\documentclass[a4paper,11pt]{article}
\usepackage[utf8]{inputenc}
\usepackage[T1]{fontenc}

\usepackage{amsthm,amsmath}
\usepackage{tikz}
\usepackage{mathrsfs,amssymb,amsfonts} 
\usepackage{enumitem}
\usepackage{fullpage}
\usepackage{hyperref, enumerate}
\usepackage[babel]{microtype}
\usepackage[english]{babel}
\usepackage[capitalise]{cleveref}

\usepackage{thmtools}
\usepackage{mathtools, comment}
\usepackage{amssymb}
\usepackage[nomath]{lmodern}
\usepackage{graphicx}
\usepackage{pgf,tikz,tkz-graph,subcaption}
\usetikzlibrary{arrows,shapes}
\usetikzlibrary{decorations.pathreplacing}
\usepackage{tkz-berge}
\usepackage{enumitem}
\usepackage[normalem]{ulem}
\usepackage{hyperref}
\hypersetup{colorlinks = true, linkcolor = blue, citecolor = blue, urlcolor = blue}

\newcommand*{\bfloor}[1]{\left\lfloor #1\right\rfloor}

\newcommand{\eps}{\varepsilon}

\allowdisplaybreaks

\usepackage[margin=1in]{geometry}
\newtheorem{defi}{Definition}
\newtheorem{conj}[defi]{Conjecture}

\newtheorem{thr}[defi]{Theorem}

\newtheorem{q}[defi]{Question}

\newtheorem{claim}[defi]{Claim}
\newcommand*{\myproofname}{Proof}
\newenvironment{claimproof}[1][\myproofname]{\begin{proof}[#1]}{\end{proof}}

\DeclareMathOperator{\lcm}{lcm}

\title{Resolution of an Erd\H{o}s' problem on least common multiples}
\author{Stijn Cambie \thanks{Department of Computer Science, KU Leuven Campus Kulak-Kortrijk, 8500 Kortrijk, Belgium. Supported by a postdoctoral fellowship by the Research Foundation Flanders (FWO) with grant number 1225224N. Email: \protect\href{mailto:stijn.cambie@hotmail.com}{\protect\nolinkurl{stijn.cambie@hotmail.com}}}}

\begin{document}
\parindent=0cm
\maketitle

\begin{abstract}
    Erd\H{o}s posed the question whether there exist infinitely many sets of consecutive numbers whose   
    least common multiple ($\lcm$) exceeds the $\lcm$ of another, larger set with greater consecutive numbers. In this paper, we answer this question affirmatively by proving that the ratio of the $\lcm$'s can be made arbitrarily large.
\end{abstract}

\hfill \small{A problem on primes presented $47$ years ago}
\normalsize

\section{Introduction}\label{sec:intro}

In his $1979$ paper (based on remarks from a conference in October 1977), Erd\H{o}s~\cite{Erdos79} posed several questions and conjectures in number theory, which he considered as unconventional.
We focus on one such simple posed question, for which the answer seemed frustratingly elusive at that time.

The least common multiple ($\lcm$) of a set of numbers is the smallest positive integer divisible by each number in the set (and thus divides their product).
The $\lcm$ of a set with larger values or greater cardinality is typically larger than the $\lcm$ of a set with smaller values or cardinality. However, this is not necessarily true for arbitrary sets.
Therefore, Erd\H{o}s considered the $\lcm$ of consecutive numbers.
For sets of two consecutive natural numbers, their $\lcm$ is larger if their values are larger, since $x<y \Rightarrow \lcm(x,x+1)=x(x+1)<y(y+1)=\lcm( y,y+1).$
In contrast, for $k \ge 3$ consecutive natural numbers, this is not true.
For example when $x$ is even and $y>x$ is odd and not much larger than $x$,
we may have that $ \lcm(x,x+1,x+2)=x(x+1)(x+2)> \frac{y(y+1)(y+2)}{2}=\lcm( y,y+1,y+2).$

The non-monotonicity can be extreme.
Let $f(n)=\lcm\{n,n+1,\ldots, n+6\}$. Then $f(1)\le f(2)\le \ldots\le f(7) \le f(8)>f(9),$ and also $f(7)=2^3\cdot3^2\cdot5\cdot7\cdot11\cdot13>f(10)=2^4\cdot3\cdot5\cdot7\cdot11\cdot13$. 
This is a counterexample to behaviour that seemed plausible to Erd\H{o}s.

When comparing the $\lcm$ of $k$ small consecutive integers with the $\lcm$ of $\ell>k $ larger consecutive integers, determining which is larger becomes harder.
For example, one might wonder which of the following two quantities is larger:
$$\lcm\{676,677,\ldots, 773\} \mbox { and } \lcm\{798,799,\ldots, 903\}.$$
Some minimal examples, with a cardinality of $7$ and smallest elements $37$ and $48$, respectively, illustrating the unexpected inequality are:
$$\lcm\{53,54,55,56,57,58,59\}>\lcm\{63,64,65,66,67,68,69,70\} \mbox{ and }$$
$$\lcm\{37,38,39,40,41,42,43,44\}>\lcm\{48,49,50,51,52,53,54,55,56\}.$$

Erd\H{o}s~\cite{Erdos79} posed the question whether there exist infinitely many cases in which the $\lcm$ of the larger such set is smaller than the $\lcm$ of the smaller set.
This question is also listed as problem $678$ at \url{https://www.erdosproblems.com/678}.
We answer the question by giving an elementary proof for the following stronger result in~\cref{sec:proof}.


\begin{thr}\label{thr:main}
    Let $C \ge 1$ be a constant.
    For every $k$ sufficiently large, there exist integers $0<x<y$ with $y>x+k$ such that
    $\lcm\{x,x+1,\ldots,x+k-1\} > C \cdot\lcm\{y,y+1,\ldots,y+k\}.$
\end{thr}

We conjecture that an even stronger form of this result is true.

 \begin{conj}\label{conj:main2}
    Let $C \ge 1$ be an arbitrary constant.
    Then there exist integers $k$ and $0<x<y$ with $y>x+k$ such that
    $\lcm\{x,x+1,\ldots,x+k-1\} > \lcm\{y,y+1,\ldots,y+k+C-1\}.$
\end{conj}

We note that a proof for~\cref{conj:main2} along the same lines as the proof in~\cref{sec:proof} would follow from density results for solutions from the Chinese remainder theorem (which even only need to hold for some values of $k$). This reduction is proven in Subsection~\ref{subsec:reduction_to_sieveQ}.

\begin{q}\label{ques:sieve-styleQ}
Let $\varepsilon>0$ be fixed.
Then there exists $0<\delta<\varepsilon$, such that for all sufficiently large $k$ the following is true.

Let the primes strictly between $\sqrt k$ and $k$ be written as $\sqrt k <p_1\le  \ldots \le p_r <k$. Let $M= \prod_{j=1}^r p_j$.
For every $1 \le j \le r$,
let $I_j= \left \{1,2,\ldots, \left \lceil p_j^{1-\delta} \right \rceil \right \}$.
Then among every $\left \lceil M^\varepsilon \right \rceil $ consecutive integers,
there is at least one value $x$ for which 
$ x \bmod p_j$ belongs to $I_j$ for every $1 \le j \le r$.
\end{q}

The state-of-the-art to solve this type of problem seems to be sieve-type bounds, but no contacted expert could tell it was possible to prove. 
The possible non-trivial nature of this (or a related) desired density result can be guessed from studies on the distribution~\cite{GK08} and equidistribution~\cite{KS21} of the outcomes from the Chinese remainder theorem.

\section{Proof of~\cref{thr:main}}\label{sec:proof}

We start with two conventions. We will assume that $k$ is taken sufficiently large such that all statements work, and will use $a \mod n$ for the unique integer between $0$ and $n-1$ which has the same residue modulo $n$ as $a.$

Let $M= \lcm\{ 1,2,\ldots,k\}=\prod_{p \le k} p^{\bfloor{ \log_p(k) }}$ be the product of all maximal prime powers bounded by $k.$ Using the notation $m=\prod_{p \le \sqrt k} p^{\bfloor{ \log_p(k) }}$, we can also write $M=m \cdot \prod_{\sqrt k <p \le k} p$.

Let $A$ and $B$ be the sets of vectors $\vec{a}= ( a_p \mid \sqrt k <p \le k )$ and $\vec{b}= ( b_p \mid \sqrt k <p \le k )$, whose coordinates satisfy $1 \le a_p \le p-(k \mod p) $ and $ p-(k \mod p) \le b_p \le p$ for every prime $k\ge p> \sqrt k$ respectively.

For every $\vec a \in A$ and $\vec b \in B$, we consider the unique solutions $0 \le x_{\vec a}, y_{\vec b} < M$ satisfying the following systems of modular equations
$$ x_{\vec a} \equiv \begin{cases}
    1 \pmod m\\
    a_p \pmod p \mbox{ for every } \sqrt k <p \le k
\end{cases} \mbox{ and } y_{\vec b} \equiv \begin{cases}
    0 \pmod m\\
    b_p \pmod p \mbox{ for every } \sqrt k <p \le k.
\end{cases}$$
The existence and uniqueness is known by the Chinese remainder Theorem.

We first prove the following claim, which says that among all linear combinations with sufficiently many choices for the coefficients, there are not too many consecutive numbers excluded.

\begin{claim}\label{clm:approx}
    Let $p_1 < p_2 <  \ldots < p_r$ be primes and $w_1, w_2, \ldots, w_r$ be integers, such that the combinations $\sum_{i} c_i w_i$ over all possible combinations with $0 < c_i \le p_i $ lead to all residues modulo $P=p_1p_2\ldots p_r.$ Let $B_i \subset [p_i]$ be a set of size at least $(1-\eps)p_i$ for every $1 \le i \le r$.
    If $\eps(p_1+p_2+\ldots+p_r)< n \le p_1,$ among every $n$ consecutive integers there is at least one which equals $\sum_{i} c_i w_i$ modulo $P$ where $c_i \in B_i$ for every $1 \le i \le r$.
\end{claim}

\begin{claimproof}
    By the condition, there is a unique solution to $\sum c_i w_i \equiv 1 \pmod P$,
    which implies for each of the $n$ consecutive integers $z$,
    $z \equiv \sum (zc_i) w_i \pmod P$.
    Among the $n$ numbers of the form $zc_i \mod p_i,$ at most $\eps p_i$ are not in the interval $B_i.$
    Thus at most $\eps(p_1+p_2+\ldots+p_r)<n$ combinations are not of the desired form, from which the conclusion trivially follows.
\end{claimproof}

By the density of primes, see e.g.~\cite{BHP01}, we know that for fixed tiny $\eps>0$, there are many primes in the range $(k/2,(1+\eps)k/2)$ and $((1-\eps)k,k)$.
We will use this to find suitable $x_{\vec a}$ and $ y_{\vec b}.$

Choose two primes $p_1, p_2 \in (k/2,(1+\eps)k/2).$
Let $B_1=[2p_1-k,p_1]=[p_1-( k \mod {p_1}),p_1]$ and $B_2=[2p_2-k,p_2]$.

Let $B'=\{ b_p=0\ \forall p \not \in \{p_1, p_2\}, b_{p_1} \in B_1, b_{p_2} \in B_2\}.$
Note that the vector $\left\{ y_{\vec{b}} \mid \vec{b} \in B'\right\} $ differ by multiples of $\frac{M}{p_1}$ and $\frac{M}{p_2}.$
By~\cref{clm:approx}, for every $(\eps(p_1+p_2))+1)\frac{M}{p_1p_2}<\frac{M}{k}$ consecutive integers, at least one can be expressed as $y_{\vec{b}}$ (for some $\vec{b} \in B')$.
This implies that there is a solution $y:=y_{\vec{b}}$ satisfying $\frac{M}{5C}(1+1/k)<y<\frac{M}{4C}-k$.

By~\cref{clm:approx} applied to $3$ primes in $((1-\eps)k,k)$, similarly as how $y$ was derived, we can choose a value $x:=x_{\vec a}$ such that $x<y<x+\frac{M}{5Ck}.$
Note that $y-x \ge m-1 > k$ by the choice of their residues modulo $m.$

\begin{claim}\label{clm:ineq_ofratios}
    $$\frac{ y\cdot(y+1)\cdot(y+2) \cdot \ldots \cdot(y+k)} { \lcm\{y,y+1,\ldots,y+k\} }= M \cdot \frac{ x \cdot (x+1)\cdot \ldots\cdot(x+k-1) }{ \lcm\{x,x+1,\ldots,x+k-1\} }$$
\end{claim}

\begin{claimproof}
    We will conclude that the left hand side equals the right hand side, $LHS=RHS,$ by proving that for every prime $p$, $v_p(LHS)=v_p(RHS)$.
    
     For any prime $p \le \sqrt k,$ we note that in $k$ or $k+1$ consecutive numbers at most one element can be a multiple of $p^{\bfloor{ \log_p(k) }+1},$ in which cases also the least common multiple of all these numbers will have the same $p$-adic value as this element.
     Taking into account that $x \equiv y+1 \equiv 1 \pmod m$, this implies that 
     \begin{align*} v_p(LHS)&= \sum_{i=y}^{y+k} \min\left( v_p(i), \bfloor{ \log_p(k) } \right) - \bfloor{ \log_p(k) } \\
     &= \sum_{i=y+1}^{y+k} \min\left( v_p(i), \bfloor{ \log_p(k) } \right)\\
     &= \sum_{i=x}^{x+k-1} \min\left( v_p(i), \bfloor{ \log_p(k) } \right)\\
     &= v_p(RHS) .
     \end{align*}
     The last inequality holds due to $ v_p(M)=\bfloor{ \log_p(k) }.$

     For every prime $\sqrt k < p \le k$, both sides contain at most one multiple of $p^2$, and so it is sufficient to count the number of multiples of $p$ in the numerators to conclude $v_p(LHS)=v_p(RHS).$
     The latter being true by the choice of $x \mod p=a_p, y \mod p=b_p$.

     Primes $p>k,$ appear at most once among $k$ or $k+1$ consecutive numbers, and so these primes satisfy $v_p(LHS)=v_p(RHS)=0.$
\end{claimproof}

Finishing the proof is now done by the following computation.

\begin{align*}
    \frac{ \lcm\{x,x+1,\ldots,x+k-1\} } { \lcm\{y,y+1,\ldots,y+k\} }&= M \cdot \frac{ x \cdot (x+1)\cdot \ldots\cdot(x+k-1) }{ y\cdot(y+1)\cdot(y+2) \cdot \ldots \cdot(y+k)  } \\ & \ge \frac{M}{y+k}  \left(\frac{x}{y}\right)^k\\
    &> 4C \left( \frac{k}{k+1}\right)^k > C.
\end{align*}

\subsection{\cref{ques:sieve-styleQ} implies~\cref{conj:main2}}\label{subsec:reduction_to_sieveQ}

Here we assume \cref{ques:sieve-styleQ} has a positive answer, and will deduce that in that case also~\cref{conj:main2} is true.

Note that by the prime number theorem, $M=k^{\Theta( k/ \log{k} )}$.
The average size of values $x$ satisfying the condition stated in \cref{ques:sieve-styleQ} among $M^{\eps}$ consecutive integers is at least $M^{\varepsilon-\delta}$, which is large for $k$ large. One can e.g. assume that $(\varepsilon-\delta) k^{0.01}> \log{k}$. The latter indicating why~\cref{ques:sieve-styleQ} should have a positive answer. 

\begin{proof}[Proof of \cref{ques:sieve-styleQ} $\Rightarrow$ \cref{conj:main2}]

Take $\eps =\frac{1}{3C+6}$ and let $0<\delta<\eps$ be a constant which works for~\cref{ques:sieve-styleQ}. Let $K$ be a sufficiently large integer.
For every prime $\sqrt K < p \le 2K$, there are at most 
$2p^{1-\delta}\left( \frac{K}{p}+1\right) \le 6 p^{-\delta} K  \le 6K^{1-\delta/2}$ values $K \le k <2K$ for which $k \mod p > p-p^{1-\delta}$ or $k \mod p < p^{1-\delta}.$
By averaging, we know that there is a choice of $K \le k <2K$ for which a fraction no more than $7K^{-\delta/2}< \eps$ of the $r$ primes $\sqrt k \le p <k$ satisfy $k \mod p > p-p^{1-\delta}$ or $k \mod p < p^{1-\delta}.$
Pick such a $k$, and let the $r$ primes be $p_1 \le p_2 \le \ldots \le p_r.$

Let $x$ and $y$ be any solution for which 
$x \mod p_j$ and $-y \mod p_j$ both belong to $I_j$ for every $1 \le j \le r$.
So except from an $\eps$ fraction of the $r$ primes these satisfy $x \mod p_j \le p_j- (k \mod p_j)$ and $y \mod {p_j} \ge p_j- (k+C-1 \mod p_j).$
For those primes, we know that 
$$v_p \left( \frac{ y\cdot(y+1)\cdot(y+2) \cdot \ldots \cdot(y+k+C-1)} { \lcm\{y,y+1,\ldots,y+k+C-1\} } \right) \ge v_p \left( \frac{ x \cdot (x+1)\cdot \ldots\cdot(x+k-1) }{ \lcm\{x,x+1,\ldots,x+k-1\} } \right)+1. $$
Let $G$ be the subset of the $r$ primes satisfying the above.
Let $B$ be the subset of primes for which it is not the case (but then it is still true with $+1$ replaced by $-1$).
Note that $\prod_{p \in B} p \le M^{2\eps}$ and $m=\prod_{p \le \sqrt k} p^{\bfloor{ \log_p(k) }} < M^{\eps}.$
As such, we have $\frac{ \prod_{p \in G} p }{m\prod_{p \in B} p }\ge M^{1-5\eps}.$
We conclude that (instead of~\cref{clm:ineq_ofratios})
$$\frac{ y\cdot(y+1)\cdot(y+2) \cdot \ldots \cdot(y+k+C-1)} { \lcm\{y,y+1,\ldots,y+k+C-1\} } >  M^{1-5\eps} \cdot \frac{ x \cdot (x+1)\cdot \ldots\cdot(x+k-1) }{ \lcm\{x,x+1,\ldots,x+k-1\} }.$$

Now by~\cref{ques:sieve-styleQ}, we can choose $x,y$ such that $M^{\eps} <y-x<3M^{\eps}$ and $M^{2\eps}<x<y+k+C<M^{3\eps}$.

The proof is now done by the following computation.

\begin{align*}
    \frac{ \lcm\{x,x+1,\ldots,x+k-1\} } { \lcm\{y,y+1,\ldots,y+k+C-1\} }&>  M^{1-5\eps} \cdot \frac{ x \cdot (x+1)\cdot \ldots\cdot(x+k-1) }{ y\cdot(y+1)\cdot(y+2) \cdot \ldots \cdot(y+k+C-1)  } \\ & > \frac{M^{1-5\eps}}{(y+k+C)^C}  \left(\frac{x}{y}\right)^k >1.
\end{align*}

\end{proof}

\section*{Acknowledgement}

The author thanks Bart Michels and Wouter Castryck for carefully reading an early draft.

\bibliographystyle{abbrv}

\end{document}